\documentclass[conference]{IEEEtran}
\IEEEoverridecommandlockouts
\usepackage{cite}
\usepackage{amsmath,amssymb,amsfonts}
\usepackage{algorithmic}
\usepackage{caption}
\usepackage{subcaption}
\usepackage{graphicx}
\usepackage{textcomp}
\usepackage{xcolor}
\usepackage{tabularx}
\usepackage{url}
\usepackage{hyperref}
\def\BibTeX{{\rm B\kern-.05em{\sc i\kern-.025em b}\kern-.08em
    T\kern-.1667em\lower.7ex\hbox{E}\kern-.125emX}}
   
\begin{document}

\title{A Transient Thermal Model for Power Electronics Systems\\
}
\author{\IEEEauthorblockN{Neelakantan Padmanabhan}
\IEEEauthorblockA{
\textit{ZF Passive Safety Systems US Inc}\\
Livonia, USA \\
neel.padmanabhan@zf.com}
}
\maketitle

\IEEEoverridecommandlockouts
\IEEEpubid{\begin{minipage}[t]{\textwidth}\ \\[10pt] 
        \centering\small{\copyright 2024 IEEE. Personal use of this material is permitted.  Permission from IEEE must be obtained for all other uses, in any current or future media, including reprinting/republishing this material for advertising or promotional purposes, creating new collective works, for resale or redistribution to servers or lists, or reuse of any copyrighted component of this work in other works.\\ 
        Accepted for publication in the IEEE SouthEastCon 24 (March 2024). \\
        DOI: 10.1109/SoutheastCon52093.2024.10500091
        }
\end{minipage}} 

\begin{abstract}
An equation based reduced order model applicable to generalized heat equation and thermal simulations of power electronics systems developed in commercial CFD tools, is presented in this work. The model considers the physics of heat transfer between multiple objects in different mediums and presents a set of equations that can be applied to a wide range of heat transfer scenarios including conduction, natural and forced convection problems. A few case studies including heat transfer in a power electronic system are simulated in Ansys\textsuperscript{\textregistered} Icepak\textsuperscript{\texttrademark} and the temperatures from the simulations are compared with the temperatures predicted by the models. The models are observed to be highly accurate when compared with the simulations. The predictive model described in this work reduces large complex simulations down to a few parameters which tremendously improves the computation speed, uses very low physical disk space and enables fast evaluation of thermal performance of the system for any changes in the input parameters. 
\end{abstract}


\section{Introduction}
Commercial CFD codes such as Ansys\textsuperscript{\textregistered} Icepak\textsuperscript{\texttrademark} are indispensable in development of simulations of electronic thermals in a wide range of industries including automotive, aerospace and energy. A high fidelity thermal simulation requires significant computational resources and time. The computational expense increases as the complexity of the system increases. Thus, for product design and development process, reduced order models (ROM) are highly desirable. Development of ROM is an active area of research \cite{ROMbook}. While Linear Time Invariant (LTI) based ROM \cite{LTI-POD} are available in Ansys\textsuperscript{\textregistered} suite of tools, they are applicable only to linear systems. Approaches for ROM development for non-linear systems include Response surface methodology (RSM) \cite{RSM}, Proper Orthogonal Decomposition (POD) \cite{POD-neel}, Subspace projection \cite{krylov} and machine learning \cite{ML} to name a few. While powerful, these approaches are computationally expensive even for simple systems. For complex systems with large number of components, the costs of simulations and subsequent ROM development increase non-linearly \cite{cost,JoT-neel}. Thus, a simple and computationally inexpensive method is advantageous. Thermal lumped element method  \cite{lumped} is a simple approach often utilized to model the thermal behavior of systems. In this method, the model constants are typically estimated from steady state behavior of a system. However, in power electronics systems, steady state assumptions often result in incorrect estimation of temperatures because the inputs are highly transient. The approach described in this work is based on the concept of thermal lumped element model, but applies modifications in how the model parameters are measured and suggests a new equation that can be used to model the transient thermal behavior of linear and non-linear systems. The ease of setup, computation speed and range of applicability make this approach significant in iterative design improvement process in development of power electronics systems. The model presented in this work is built in steps and presented in the following sections. A model for a single body insulated system is built from the transient form of heat conduction equation. The model is then extended to two body systems and the appropriate modifications are presented. This is then extended to a multi-body system with mixed modes of heat transfer. A method to estimate the model constants for a general system and evaluation of the model with simulation data are then presented. 
\section{Background}
The generalized heat conduction equation can be expressed as, 
\begin{equation}
    \rho c_p \frac{\partial T}{\partial t} = \frac{\partial}{\partial x_i} (k \frac{\partial T}{\partial x_i}) + P,
\end{equation}
where, $P$ is the source (input) term, $\rho,c_p,k$ represent the density, specific heat, and thermal conductivity, respectively, $T$ the temperature, and $x_i,t$ the spatial and temporal terms. For constant properties at atmospheric pressures \cite{constprops}, this equation can be expressed as,
\begin{equation}
    \frac{\partial T}{\partial t} = \alpha \nabla^2 T + \frac{P}{\rho c_p},
\end{equation}
where, $\alpha=k/\rho c_p$, is the thermal diffusivity and $\nabla^2$ is the Laplace operator. Depending on the type of heat transfer across the surface of the solid bodies, the boundary conditions can be formulated as insulated walls with no heat flux, conduction across different materials or convection between a solid surface and the surrounding fluid. Integrating this equation with respect to time results in an equation of the form,
\begin{equation} \label{eqmain}
    T(t)=D(t) + \frac{P}{\rho c_p} t + c_{int}. 
\end{equation}
where, $D(t)=\alpha \int \nabla^2 T dt$ is the diffusion term and $c_{int}$ is the integration constant. This form of the equation is used as a starting point to develop models for the thermal systems.
\section{Models}
\subsection{Model for 1-body insulated system}
For one body insulated system, it is observed that the diffusion of temperature across the body is instantaneous. Thus, the diffusion term (spatial gradient) in Eq. \ref{eqmain} is negligible. Hence, evolution of temperature of the body can be modeled as,
\begin{equation}
    T^B(t)=T^0 + T^L(t) = T^0 + \Bigl(\frac{P^B}{C^B} \Bigr) t,
\end{equation}
where, $T^0$ is the initial temperature, $C^B=\rho c_p V$ is the thermal capacitance of the body, $V$ the volume of the body and $P^B$ the power input applied to the body. This temperature curve is observed to be linear with the slope $P^B/C^B$. Thermal simulations are created for 1-body insulated system where wall boundary conditions with zero heat flux are applied at the cabinet. The temperature from the simulations and the temperature from the model are presented in Figure \ref{1-body_comp}.
\subsection{Model for 2-body conduction heat transfer}
\begin{figure}
     \centering
     \begin{subfigure}[b]{0.3\textwidth}
         \centering
         \includegraphics[width=\textwidth]{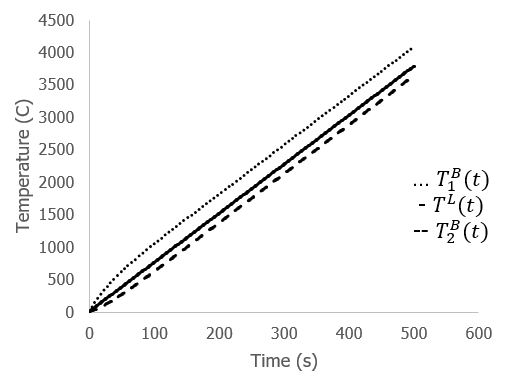}
     \end{subfigure}
     \hfill
     \begin{subfigure}[b]{0.3\textwidth}
         \centering
         \includegraphics[width=\textwidth]{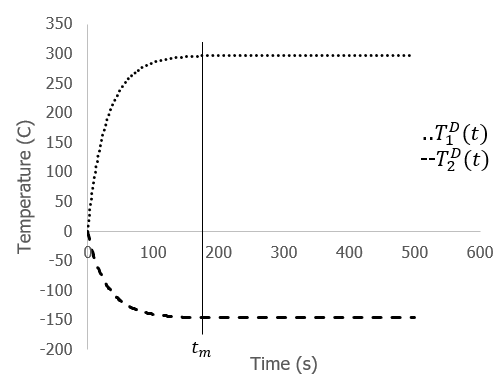}
     \end{subfigure}
        \caption{Temperature of bodies, $T^L(t)$ from 2-body simulation and deviation terms}
        \label{TD_model}
\end{figure}

\begin{figure}
\centerline{\includegraphics[width=0.5\textwidth]{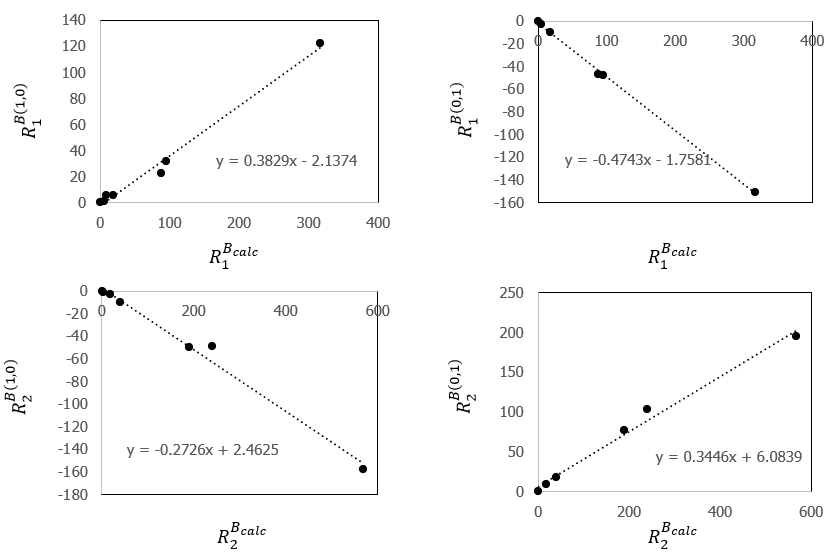}}
\caption{Estimated vs calculated thermal resistance}
\label{Rcalc_linear}
\end{figure}

For heat conduction between 2 bodies insulated from the surrounding, it is observed that there is a spatial gradient in temperature from the heat source to the sink. Further, the heat transfer in this case can be approximated to 1D heat transfer since the temperature gradient is observed to be significant only in the direction of heat flow. The temperatures (at center) of each body, $T^B _{1} (t)$ and $T^B_{2} (t)$, can thus be modeled as,
\begin{equation}
\begin{aligned}
    T^B _{1} (t)=T^0+T^L(t)+T^{D_{m}} _{1}(t),\\ 
    T^B _{2} (t)=T^0+T^L(t)+T^{D_{m}} _{2}(t),
\end{aligned}
\end{equation}
where, $T_0$ denotes the initial temperature, $T^L(t)=(P^{B} _{T}/C^{B} _{T})t$ and the term $P^{B} _{T}/C^{B} _{T}$ denotes the ratio of total applied power in both bodies to the total thermal capacitance of both bodies. For two body system it is defined as,
\begin{equation}
    \frac{P^{B} _{T}}{C^{B} _{T}} = \frac{\sum_{i=1} ^2 P^{B} _{i}}{\sum_{i=1} ^2 C^{B} _{i}}.
\end{equation}
The significance of the term $T^L (t)$ is that, it represents a temperature curve that lies between $T^{B} _{1}(t)$ and $T^{B} _{2}(t)$ and is always linear with a slope of $P^{B} _{T}/C^{B} _{T}$. This term can be computed a priori using the material properties and the input powers. The location $x_{P^{B} _{T}/C^{B} _{T}}$, within the system at which this temperature curve is observed is variable and depends on the applied power, geometry, and material properties. For similarly sized bodies with similar material properties, $x_{P^{B} _{T}/C^{B} _{T}}$ is typically located close to the midpoint between the two bodies. However, when there are differences in thermal capacitance of bodies, this location is closer to body with the larger thermal capacitance. 
The second term, $T^{D_{m}} _{i}(t)$ is a model for $T^{D} _{i}(t)$ which is the deviation in temperature of the body measured from the simulation $T^{{B_{sim}}} _{i} (t)$ and the linear temperature curve $T^L(t)$. This is defined as,
\begin{equation}
    T^{D} _{i}(t)=T^{{B_{sim}}} _{i}(t)-T^L(t),
\end{equation}
It is observed that $T^{D} _{i}(t)$ (Fig. \ref{TD_model}) has the shape of a characteristic exponential curve and can be modeled similar to transient response of an RC circuit \cite{transRC1,transRC2,transRC3}. The model equation can be expressed as, 

\begin{equation}
    T^{D_{m}} _{i}(t)=A_i (1-e^{-k_i t}),
\end{equation}
where, $A_i=|T^{B_{sim}} _{i}(t_{m})-T^L(t_{m}))|$ is the magnitude of difference in stationary temperature of each body and the linear temperature curve, and $t_{m}$ the time instant after which the temperature curve has reached a stationary state. The sign of $A_i$ is positive if  $T^{B_{sim}} _{i} > T^L$ and vice-versa. Further, the term $A_i$ can be modeled as,
\begin{equation}
    A_i = P^{B} _{i} (x) R^{B} _{i},
\end{equation}
where, $P^{B} _{i}(x)$ is the applied input power at each body. The input power is described as a spatial function because the applied power varies spatially across the system. $R^{B} _{i}$ is the thermal resistance between the body’s center and the location $x_{P^{B} _{T}/C^{B} _{T}}$. The term $k_i=1/(C^{B} _{T} R^{B} _{i})$ represents the time constant or the rate of change of the temperature. 
\begin{figure}
\centerline{\includegraphics[width=0.5\textwidth]{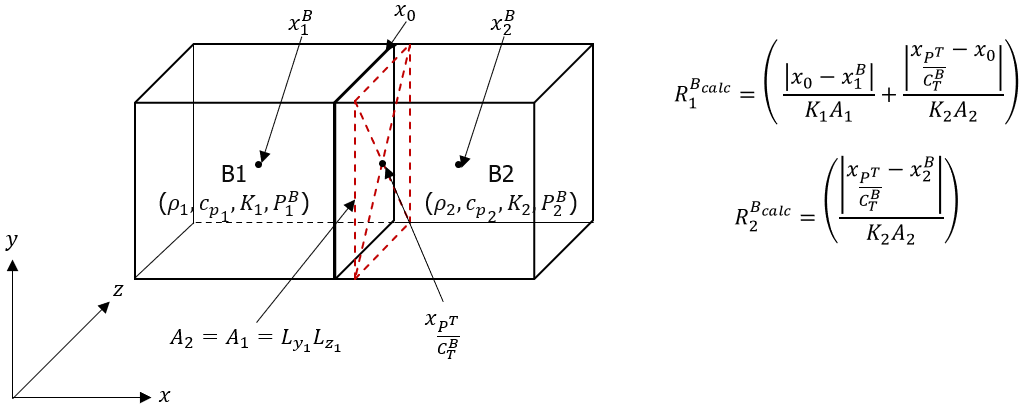}}
\caption{Analytical computation of thermal resistance}
\label{R_calc}
\end{figure}
To model $T^{D_{m}} _{i} (t)$ for a general 2-body system as described in Fig. \ref{R_calc}, where the source and sink may vary, the principle of superposition is applied. It is expressed as,
\begin{equation}
    \begin{aligned}
        T^{D_{m}} _{1}(t) = (P^{B}_{1} R^{B{(1,0)}} _{1} + P^{B} _{2} R^{B{(0,1)}} _{1}) [1-e^{-k_1 t}], \\
        T^{D_{m}} _{2}(t) = (P^{B}_{1} R^{B{(1,0)}} _{2} + P^{B} _{2} R^{B{(0,1)}} _{2}) [1-e^{-k_2 t}],
    \end{aligned}
\end{equation}
where, $R^{B{(1,0)}} _{1}$ is the characteristic thermal resistance between the center point of body $B_1 (x^{B} _{1})$ and the location $x_{P^{B} _{T}/C^{B} _{T}}$, when input $P^{B} _{1}=1$ is applied to body $B_1$ and input $P^{B} _{2}=0$ is applied to body $B_2$. This term can be determined by running a parametric simulation as described in Sec. \ref{sec_multibody} and determining the difference in the temperature between the body and $T^L$, 
\begin{equation}
R^{B{(1,0)}} _{1}=\frac{T^{B_{sim} {(1,0)}} _{1} (t_{m}) - T^L(t_{m})}{P^{B} _{1}}.
\end{equation}
In this equation, $T^{B_{sim} {(1,0)}} _{1} (t_{m})$ is the temperature of $B_1$ when input $P^{B} _{1}=1$ is applied to $B_1$ and input $P^{B} _{2}=0$ is applied to $B_2$. Since a unit input power is applied, the magnitude of temperature difference is equivalent to the thermal resistance. Similarly,
\begin{equation}
    \begin{aligned}
        R^{B{(0,1)}} _{1} = T^{B_{sim} {(0,1)}} _{1} (t_{m}) - T^L (t_{m}); P^{B} _{1}=0, P^{B} _{2}=1, \\
        R^{B{(1,0)}} _{2} = T^{B_{sim} {(1,0)}} _{2} (t_{m}) - T^L (t_{m}); P^{B} _{1}=1, P^{B} _{2}=0, \\
        R^{B{(0,1)}} _{2} = T^{B_{sim} {(0,1)}} _{2} (t_{m}) - T^L (t_{m}); P^{B} _{1}=0, P^{B} _{2}=1.
    \end{aligned}
\end{equation}
The time constants are expressed as, $k_1=1/(R^{B{(1,0)}} _{1} C^{B} _{T})$, and $k_2=1/(R^{B{(0,1)}} _{2} C^{B} _{T})$. The thermal resistance terms can also be determined analytically for simple configurations. A general form of thermal resistance for conduction is expressed as,
\begin{equation} 
    R^{B_{calc}} _{i} = \frac{L_i}{K_i A_i} = \frac{x^{B}_{i} - x_{P^{B} _{T} /C^{B} _{T}}}{K_i A_i},
\end{equation}
where, the location $x_{P^{B} _{T} /C^{B} _{T}}$ can be approximated close to the center of the two bodies. An illustration of this is provided in Fig. \ref{R_calc}. A set of trials are conducted by varying the material properties, geometry, and sizes of each body. Further the sources and sinks of heat for each case are also varied. It is observed that a linear relationship exists between the thermal resistance values determined from optimization and thermal resistance calculated analytically, ($R^{B} _{i} = m_i R^{B_{calc}} _{i} + c_i$) as observed in Fig. \ref{Rcalc_linear}.

Further, the changes in thermal resistance terms $R^{B{(1,0)}} _{1}, R^{B{(1,0)}} _{2}, R^{B{(0,1)}} _{1}, R^{B{(0,1)}} _{2}$ due to changes in geometry or material can be calculated by scaling the terms $L,k,A$. For two body system with transient powers, $P^{B} _{1}(t), P^{B} _{2}(t)$, a piece-wise approximation method is applied. The model equations are expressed as follows, 
\begin{equation}
    T^L (t)= \sum_j ^{n_T} \left[ \frac{P^{B} _{T,j}}{C^{B} _{T}} - \frac{P^{B} _{T,j-1}}{C^{B} _{T}} \right] (t-t^{0} _{j-1}),
\end{equation}
with $P^{B} _{T,j}=P^{B} _{1,j}+P^{B} _{2,j}$. 
\begin{equation}
\begin{split}
    T^{D_{m}} _{1}(t) = \sum_j ^{n_T} [(P^{B} _{1,j} R^{B{(1,0)}} _{1} + P^{B} _{2,j} R^{B{(0,1)}} _{1}) \\ - (P^{B} _{1,j-1} R^{B{(1,0)}} _{1} + P^{B} _{2,j-1} R^{B{(0,1)}} _{1})] \\ (1-e^{-k_1 (t-t^{0} _{j-1})}), \\
    T^{D_{m}} _{2}(t) = \sum_j ^{n_T} [(P^{B} _{1,j} R^{B{(1,0)}} _{2} + P^{B} _{2,j} R^{B{(0,1)}} _{2}) \\ - (P^{B} _{1,j-1} R^{B{(1,0)}} _{2} + P^{B} _{2,j-1} R^{B{(0,1)}} _{2})] \\ (1-e^{-k_2 (t-t^{0} _{j-1})}),
\end{split}
\end{equation}
where, $j$ is the index of number of changes in the input power, $n_T$ is total number of changes (transients) in the input, and $t^{0} _{j}$ is the time instant at which the change in input powers occurs. The initial values are, $P^{B} _{1,0}=P^{B} _{2,0}=0, t^{0} _{0}=0$. Thermal simulations are created for 2-body insulated system where wall boundary conditions with zero heat flux are applied at the cabinet. The temperature from the simulations and the temperature from the model are presented in Fig \ref{2-body_comp}.
\subsection{Model for convection heat transfer}
For convection heat transfer, with one solid body, the problem essentially has two mediums, the solid and the surrounding fluid. When the flow velocity of the fluid is non-zero it is observed that the temperature evolution of the solid body follows an exponential curve. Thus, for one solid body in a fluid, the temperature evolution of the solid body is given by,
\begin{equation}
    T^{B-F}(t)=T^0 + P^{B} _{1} R^{B-F} (1-e^{\frac{t}{R^{B-F} C^{B} _{T}}}).
\end{equation}
The term $R^{B-F}=1/hA$ denotes the thermal resistance between the body and the fluid, $h$ the convection heat transfer coefficient, $A$ the cross-sectional area, and $C^{B} _{T}=C^{B} _{S}+C^{B} _{F}$ is the total thermal capacitance of the solid body and the fluid. While determination of convection heat transfer coefficient for complex geometries typically requires multiple experimental runs, it can also be calculated from Ansys\textsuperscript{\textregistered} Icepak\textsuperscript{\texttrademark} for a given configuration \textit{[Report $>$ Full Report {Variable: Convection Heat Transfer Coefficient}]}. A set of trial runs were performed at different flow velocities, where the $R^{B-F}$ that best fits the characteristic temperature $T^{B-F} (t)$ for the given configuration, were estimated using curve-fitting. Further, the convection heat transfer coefficient and area calculated from Icepak\textsuperscript{\texttrademark} for the given configuration were determined. The measured and estimated values of heat transfer coefficients are compared in Table \ref{tab:h}. It is observed that the heat transfer coefficient estimated from curve-fit is same as the heat transfer coefficient obtained from the tool. This provides a confirmation of the form of the model equation that best captures the heat transfer in this configuration. For simple geometries and flow configurations the convection heat transfer coefficient can also be estimated analytically \cite{HMT,HMT2}.
\begin{table}[h]
\begin{center}
\caption{$R^{B-F}$ represents the thermal resistance estimated from curve-fitting $T^{B-F}(t)$, $A$ the surface area measured from the configuration, and $h^{calc}$ the convection heat transfer coefficient calculated by Ansys\textsuperscript{\textregistered} Icepak\textsuperscript{\texttrademark}.} 
\begin{tabular}{|c|c|c|c|c|}
\hline
\begin{tabular}[c]{@{}c@{}}$U^{flow}$\\ $(m/s)$\end{tabular} & \begin{tabular}[c]{@{}c@{}}$R^{B-F}=$\\ $1/h^{est}A$\\ $(K/W)$\end{tabular} & $A (m^2)$ & \begin{tabular}[c]{@{}c@{}}$h_{est}=$\\ $1/A R^{B-F}$\\ $(W/m^2 K)$\end{tabular} & \begin{tabular}[c]{@{}c@{}}$h^{calc}$\\ $(W/m^2 K)$\end{tabular} \\ \hline
0.00098                                                      & 303.692                                                                     & 0.00015   & 21.95                                                                            & 21.83                                                            \\ \hline
5                                                            & 118.9                                                                       & 0.00015   & 56.07                                                                            & 56.19                                                            \\ \hline
50                                                           & 28.47                                                                       & 0.00015   & 234.09                                                                           & 235.7                                                            \\ \hline
\end{tabular} 
\label{tab:h}
\end{center}
\end{table}
\subsection{Model for multiple bodies with mixed conduction and convection heat transfer} \label{sec_multibody}
For a general multi-body system with mixed conduction and convection heat transfer modes, the form of model equations remains consistent with the two-body conduction / convection models. The equations for the system are expressed as,
\begin{equation}
    \begin{aligned}
    T^L (t) = \frac{P^{B} _{T}}{C^{B} _{T}} t = \frac{\sum_{i} ^{N_B} P^{B} _{i}}{\sum_{i} ^{N_B} C^{B} _{i}} t, \\
    T^{B} _{i} (t) = T^0 + T^L (t) + T^{D_{m}} _{i} (t), \\
    T^{D_{m}} _{i}(t) = \sum_i ^{N_B} P^{B} _i R^{B_{{char}}} _i [1-e^{-t/R^{B_{{char}}} _i C^{B} _{T}}].
    \end{aligned}
\end{equation}
The term $N_B$ is the total number of bodies in the system, including the fluid, $R^{B_{{char}}} _i$ the characteristic thermal resistance between the body $i$ and location $x_{P^{B} _{T}/C^{B} _{T}}$. Further, the characterization of thermal resistances are required only for the bodies / sub-systems of interest. The temperature from the simulations and the temperature from the model are presented in Fig \ref{ECU}
\subsubsection{Estimation of thermal resistance and time constants for general multibody system}
\begin{itemize}
    \item A parametric study is performed by applying 1W of input power to each source, while keeping the inputs at other bodies at 0. The number of trials is equal to the number of sources. The temperature response of each body is measured for each trial. Each simulation is run only for a short duration ($t_m \approx 20 s$), since the transient inputs in power electronic systems typically change every few seconds. The model parameters that characterize the temperature curve in the shorter time range result in higher accuracy. The parametric study can be automated in Ansys\textsuperscript{\textregistered} Icepak\textsuperscript{\texttrademark} from the \textit{Parameters and Optimization $>$ Parametric trial} option.
    \item The slope $P^{B} _{T}/C^{B} _{T}$  and the linear temperature curve $T^L (t)$ are computed. For complex geometries the volume of the body can be determined using a mechanical CAD package or by approximating the complex geometry with a number of smaller cuboid geometry whose volume can be determined easily.
    \item The deviations $T^{D} _{i}(t)= T^{B_{sim}} _{i} (t)-T^L (t)$, in simulation temperatures from each trial and the linear temperature are measured.
    \item The characteristic thermal resistance ($R^{B_{char}} _i=A_i$) and time constants ($k_i$) for each body are determined,
    \begin{itemize}
        \item The parameters, $A_i$ and $k_i$ that minimize the sum of square error between the temperature deviation $T^{D} _{i}(t)$ and the model temperature $T^{D_{m}} _{i}(t)$, $SSE=\sum(T^{D} _{i}-T^{D_{m}} _{i})^2$, are estimated using a non-linear optimization algorithm. In this work a non-linear Generalized Reduced Gradient solver in Excel\textsuperscript{\texttrademark} was used to evaluate the parameters. Estimation was also performed with the \textit{fminsearch} function in Matlab \textsuperscript{\texttrademark} and similar parameters were obtained. The sign of $A_i$ is negative if $T^{D} _{i}-T^L<0$.
    \end{itemize}
    \item The model deviation terms $T^{D_{m}} _{i} (t)$ are computed and the model temperature of the bodies are then computed as, $T^{B} _{i}(t)=T^0+T^L(t)+T^{D_{m}} _{i} (t)$.
\end{itemize}
\section{Simulation Setup}
An overview of the simulation setup used in this study to validate the models can be found in this work \cite{neel-ecce}. The simulations include one body, two body and multi-body systems in insulated, pure conduction, natural convection and forced convection modes. Detailed mechanical CAD models, electrical CAD models and material properties are used for the components (PCBA, MOSFET, Shunt) in simulation of power electronics system. All the simulations are run at pressure of $1 atm$ and ambient temperature of $20 ^{\circ}C$. For natural convection problems, the Boussinesq approximation is applied. The grid resolution in every coordinate direction is set such that the measured statistics become independent of the resolution.
\section{Results}
\begin{figure}
     \centering
     \begin{subfigure}[b]{0.24\textwidth}
         \centering
         \includegraphics[width=\textwidth]{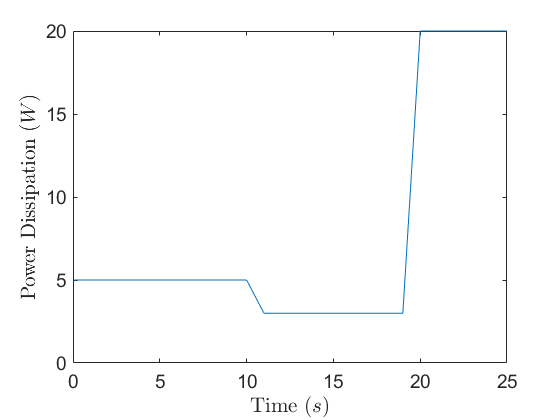}
         \subcaption{Input power dissipation}
     \end{subfigure}
     \begin{subfigure}[b]{0.24\textwidth}
         \centering
         \includegraphics[width=\textwidth]{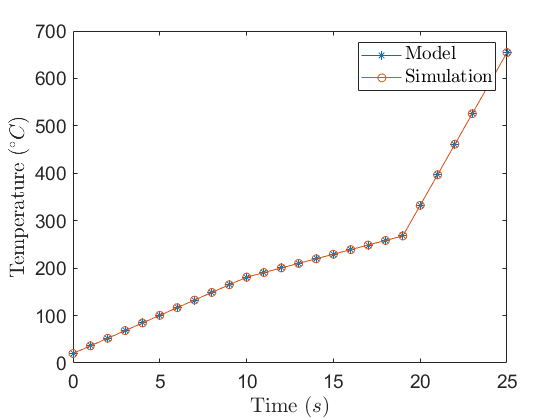}
        \subcaption{Material: Silver (metal)}
     \end{subfigure}
     \hfill
     \begin{subfigure}[b]{0.24\textwidth}
         \centering
         \includegraphics[width=\textwidth]{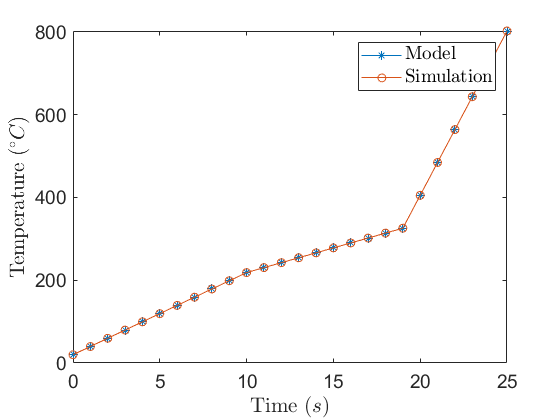}
         \subcaption{Material: Asbestos (insulator)}
     \end{subfigure}
     \begin{subfigure}[b]{0.24\textwidth}
         \centering
         \includegraphics[width=\textwidth]{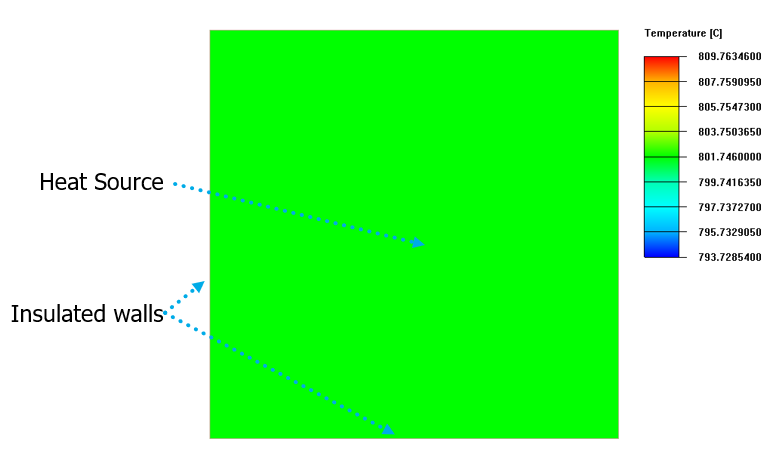}
         \subcaption{Simulation configuration - 1-body (Asbestos)}
     \end{subfigure}     
     \caption{Model vs simulation for 1-body insulated system}
\label{1-body_comp}        
\end{figure}

\begin{figure}
     \centering
     \begin{subfigure}[b]{0.24\textwidth}
         \centering
         \includegraphics[width=\textwidth]{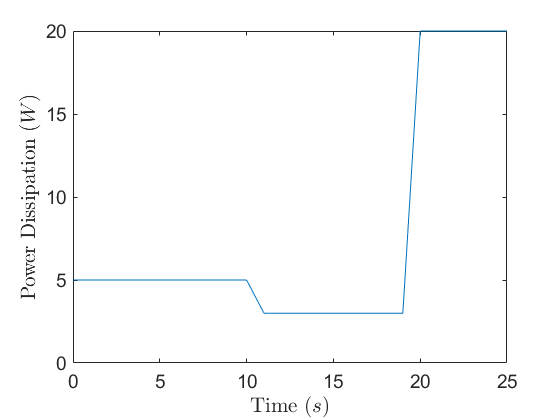}
         \subcaption{Input power dissipation}
     \end{subfigure}
     \begin{subfigure}[b]{0.24\textwidth}
         \centering
         \includegraphics[width=\textwidth]{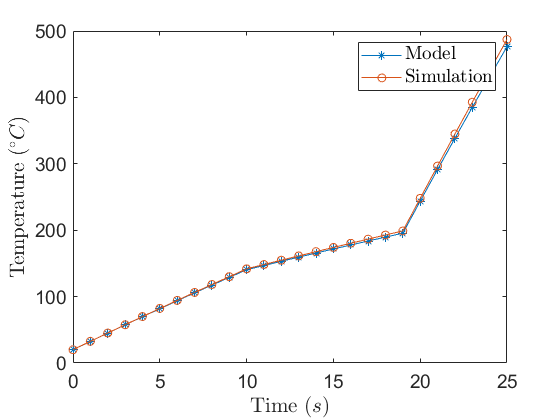}
        \subcaption{Natural convection}
     \end{subfigure}
    \begin{subfigure}[b]{0.24\textwidth}
         \centering
         \includegraphics[width=\textwidth]{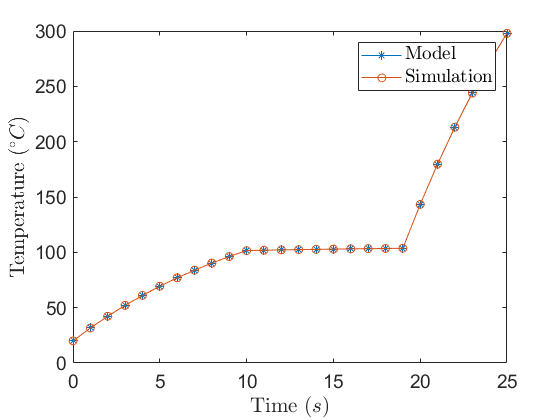}
         \subcaption{$U_{flow}=50 m/s$}
     \end{subfigure}
     \begin{subfigure}[b]{0.24\textwidth}
         \centering
         \includegraphics[width=\textwidth]{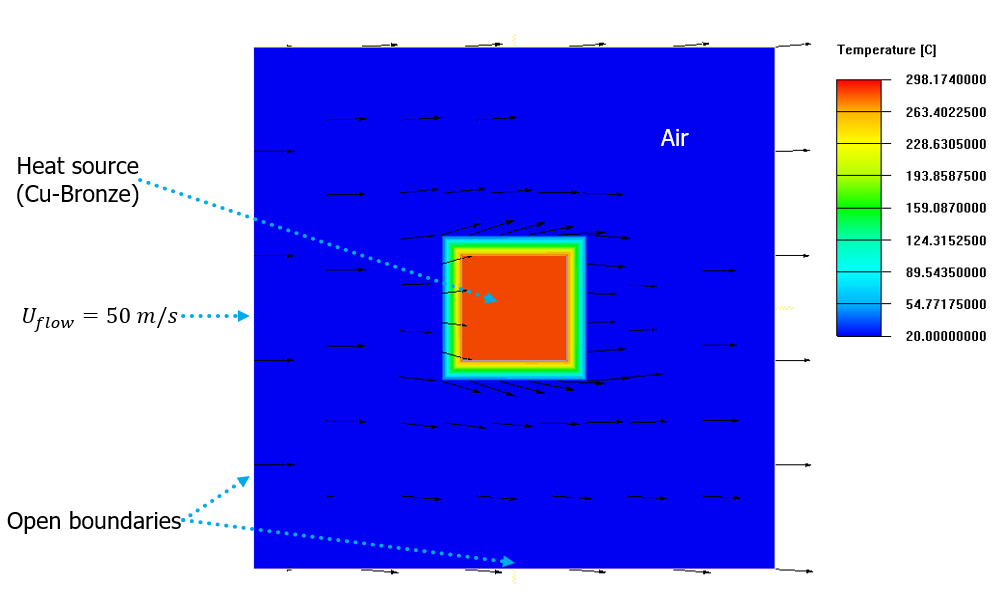}
         \subcaption{Simulation configuration - Forced convection}
     \end{subfigure}     
     \caption{Model vs simulation for 1-solid body (Cu-Bronze), fluid (air) medium natural and forced convection system}
\label{conv}        
\end{figure}

\begin{figure}
     \centering
     \begin{subfigure}[b]{0.24\textwidth}
         \centering
         \includegraphics[width=\textwidth]{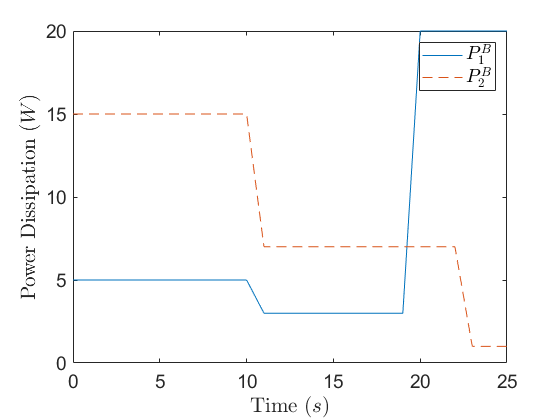}
         \subcaption{Input power dissipation}
     \end{subfigure}
     \begin{subfigure}[b]{0.24\textwidth}
         \centering
         \includegraphics[width=\textwidth]{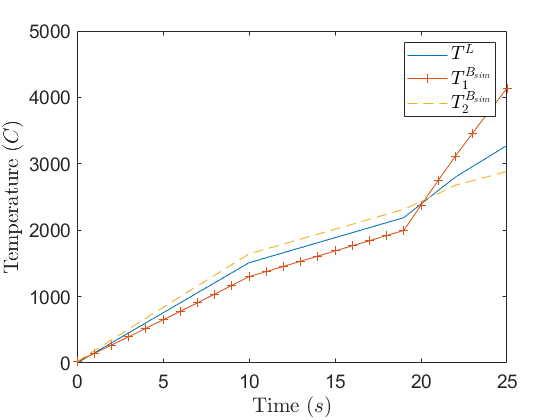}
        \subcaption{$T^L(t)$ computed a priori}
     \end{subfigure}
     \hfill
     \begin{subfigure}[b]{0.24\textwidth}
         \centering
         \includegraphics[width=\textwidth]{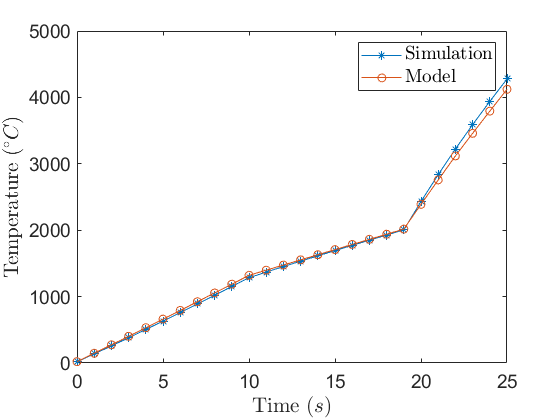}
         \subcaption{Temperature of body 1}
     \end{subfigure}
    \begin{subfigure}[b]{0.24\textwidth}
         \centering
         \includegraphics[width=\textwidth]{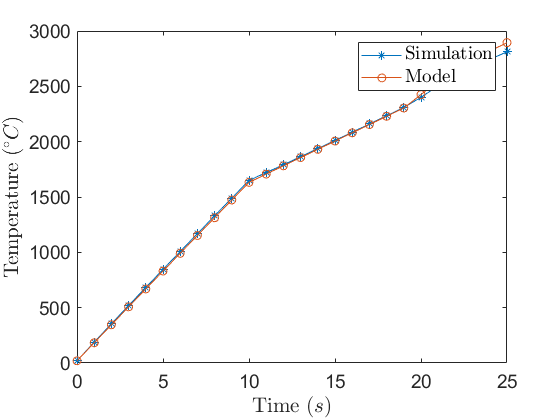}
         \subcaption{Temperature of body 2}
     \end{subfigure}
    \begin{subfigure}[b]{0.24\textwidth}
         \centering
         \includegraphics[width=\textwidth]{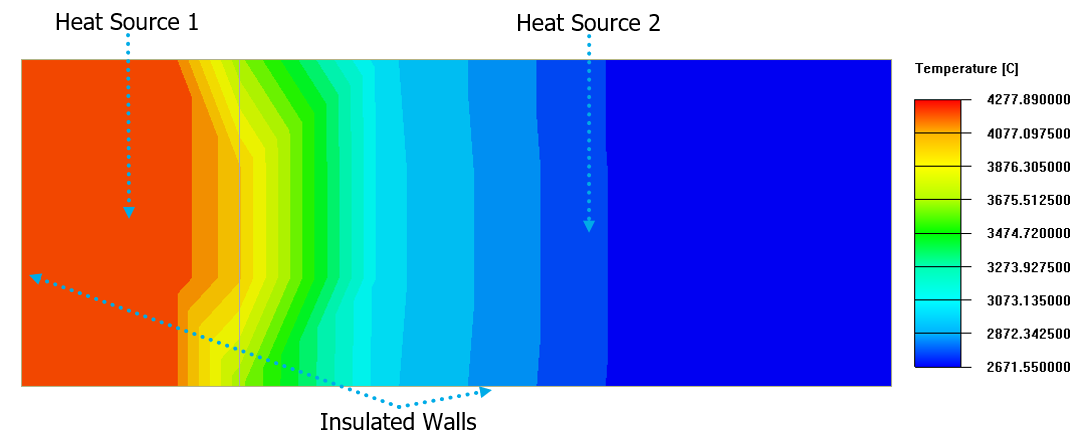}
         \subcaption{Simulation configuration}
     \end{subfigure}     
     \caption{Model vs simulation for 2-body conduction. Material 1: Ag, Material 2: FR4}
\label{2-body_comp}        
\end{figure}

\begin{figure}
     \centering
     \begin{subfigure}[b]{0.24\textwidth}
         \centering
         \includegraphics[width=\textwidth]{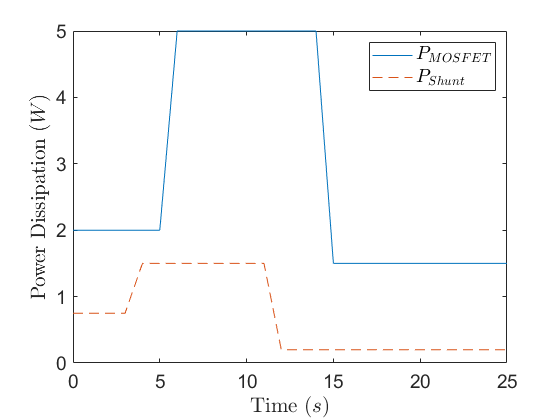}
         \subcaption{Input power dissipation}
     \end{subfigure}
     \begin{subfigure}[b]{0.24\textwidth}
         \centering
         \includegraphics[width=\textwidth]{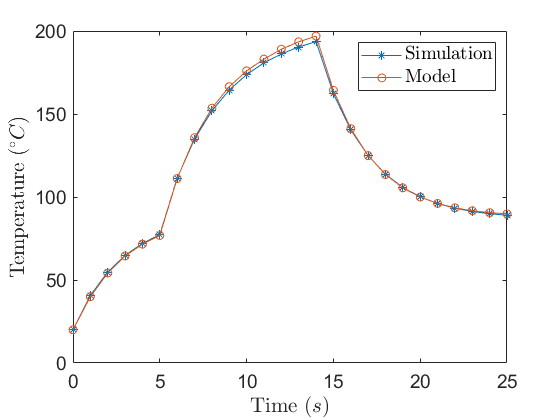}
        \subcaption{Temperature of MOSFET}
     \end{subfigure}
     \hfill
     \begin{subfigure}[b]{0.24\textwidth}
         \centering
         \includegraphics[width=\textwidth]{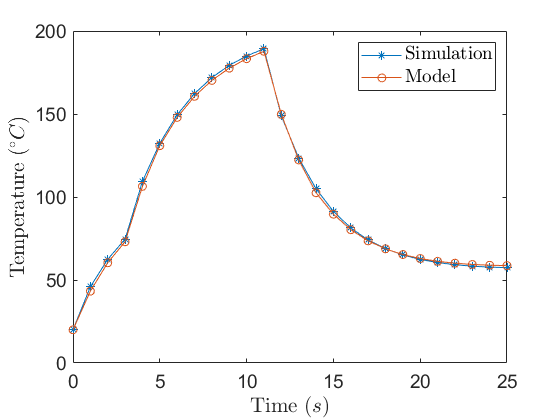}
         \subcaption{Temperature of shunt}
     \end{subfigure}
    \begin{subfigure}[b]{0.24\textwidth}
         \centering
         \includegraphics[width=\textwidth]{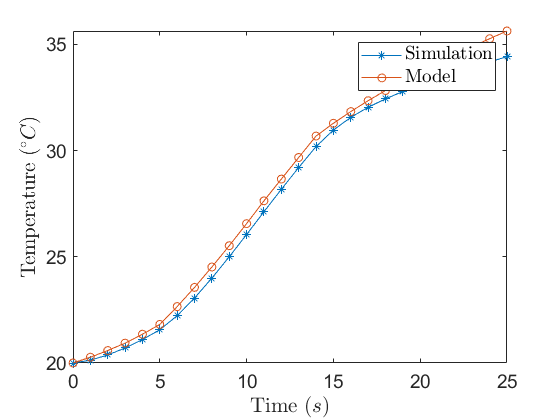}
         \subcaption{Temperature of multilayer PCBA}
     \end{subfigure}
    \begin{subfigure}[b]{0.3\textwidth}
         \centering
         \includegraphics[width=\textwidth]{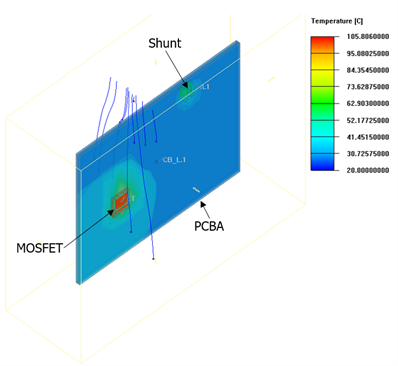}
         \subcaption{Simulation configuration}
     \end{subfigure}     
     \caption{Model vs simulation for a general multi-body power electronics system - PCBA with MOSFET and Shunt in a natural convection environment}
\label{ECU}        
\end{figure}

Comparison between the model and simulation for 1-body with insulated boundaries is presented in Fig. \ref{1-body_comp}, comparison for convection is presented in Fig. \ref{conv}, comparison for 2-body conduction is presented in Fig. \ref{2-body_comp} and the comparison for a general multi-body configuration (for power electronics) with mixed conduction and convection modes is presented in Fig. \ref{ECU}. For brevity, the results from the finite difference codes for generalized heat equation are not presented in this paper. The models are observed to be highly accurate for each case with the calculated average error of less than $3\%$. A small deviation in the simulation and model temperature of the PCBA in Fig. \ref{ECU} is observed. This is seen because the PCBA is a body with high thermal capacitance and the temperature deviation measured between the PCBA and the linear temperature curve for the simulation time range is not purely exponential. The model files occupy $0.01 \%$ of the total physical disk space that detailed simulation and solution files typically occupy. 
\begin{table}[h] 
\begin{center}
\caption{Comparison of run time between simulations in Icepak \textsuperscript{\texttrademark} and model in Matlab \textsuperscript{\texttrademark} (simulations and model code were run on the same computer). It is to be noted that parametric simulations are required to obtain the model coefficients before running the model. The parametric simulations however are run at a constant (large) $\Delta t$ and are required to be run only once.}
\begin{tabular}{|c|c|c|c|}
\hline
Case                           & \begin{tabular}[c]{@{}c@{}}Simulation \\ run time \\ in Icepak\\ (s)\end{tabular} & \begin{tabular}[c]{@{}c@{}}Parametric\\ Simulation \\ run time \\ in Icepak \\ (s)\end{tabular} & \begin{tabular}[c]{@{}c@{}}Model \\ run time \\ in Matlab\\ (s)\end{tabular} \\ \hline
1-Body Insulated               & 130                                                                               & 0                                                                                               & 0.072679                                                                     \\ \hline
1-Body Convection              & 137                                                                               & 120                                                                                             & 0.22                                                                         \\ \hline
2-Body Conduction              & 122                                                                               & 370                                                                                             & 0.1790                                                                       \\ \hline
Power electronics simulation & 126                                                                               & 240                                                                                             & 0.229                                                                        \\ \hline
\end{tabular}

\label{tab:comptime}
\end{center}
\end{table}

\section{Model limitations and correction for steady state systems}
The approach presented in this work, models the deviation $T^{D}(t)$ in temperature between a body and the linear temperature curve. This approach is however applicable only to transient systems where the inputs vary rapidly. For steady state systems where a body's temperature $T^B(t)$ follows an exponential curve, temperature deviation cannot be modeled in the same manner as described in this work. For such a system, $T^B$ can be directly modeled as $T^B(t)=T^0 + \sum_i ^{N_B} P^B _i R^{char} _i [1-e^{-k_i t}]$. 
\section{Conclusions}
The approach described in this work can accurately model the temperatures from simulations of highly transient power electronics systems in various heat transfer modes. Future work based on this approach includes evaluation of applicability of the model to experimental test setups of power electronics systems with mixed heat transfer modes, and including the effects of radiation. This includes development of analytical models for thermal resistance and capacitance of various sub-components, measurement of thermal parameters from test setups and validation of model temperatures against temperatures from real systems.

\appendix
A code of the model in Matlab \textsuperscript{\texttrademark} for various configurations is presented in the following link: \url{https://github.com/neelp-87/ROM-Transient-thermal-model-for-Power-Electronics-}


\vspace{12pt}
\color{red}

\end{document}